\documentclass[leqno,a4paper,twoside,11pt]{article}

\usepackage[T1]{fontenc}
\usepackage[latin1]{inputenc}
\usepackage[french]{babel}
\usepackage[active]{srcltx}

\usepackage{bbm,pifont,eufrak,eucal,mathrsfs,latexsym,graphicx,epsfig,color,rotating,multirow}

\usepackage{hyperref}

\usepackage[all]{xy}
\xyoption{2cell}

\newcounter{infra}[page]

\newenvironment{dem}[1][]{%
{\bf D\'emonstration #1 : }}{%
\hspace*{\fill}\nolinebreak[1]\hspace*{\fill}\underline{\bf Q.e.d.}\\}

\newenvironment{dem*}[1][]{%
{\bf D\'emonstration #1 : }}{%
 }

\newenvironment{eq*}{\begin{eqnarray*}}{\end{eqnarray*}}

\newtheorem{thm}{Th\'eor\`eme}[section] \newtheorem{lem}[thm]{Lemme}
\newtheorem{pro}[thm]{Proposition} \newtheorem{cor}{Corollaire}[thm]

\newcommand{\croi}{\times}

\newcommand{\adh}{\overline}

 \newcommand{\boufl}{\xymatrix{ \!\!
\ar@(ur,dr) }}

\newcommand{\cad}{\mbox{\it c-\`a-d }}
\newcommand{\call}{\mathscr}

\newcommand{\cf}{{\it cf. }}

\newcommand{\ch}{\vee}

\newcommand{\cten}{\Box
\hspace*{-1.77ex}\raisebox{.19ex}{$\croi$}}

\newcommand{\diag}{\mathrm{diag}}

\newcommand{\donne}{\mapsto}

\newcommand{\equi}{\Leftrightarrow} \newcommand{\et}{\mathrm{\; et \;}}

   \newcommand{\goth}{\mathfrak}

 \newcommand{\ie}{\emph{i.e. }}
\newcommand{\im}{\mathrm{Im\:}} \newcommand{\impliq}{\Rightarrow}

\newcommand{\incl}[1][r]
  {\ar@<-0.2pc>@{^(-}[#1] \ar@<+0.2pc>@{-}[#1]}

\renewcommand{\int}{\mathrm{int\:}}

 \newcommand{\inv}{^{-1}}
 \newcommand{\iso}{\simeq}

\newcommand{\kk}{{\mathbf{k}}}

\newcommand{\ma}{\displaystyle} 

\newcommand{\moins}{\:\setminus\:} 

\newcommand{\y}{y}

 \newcommand{\Ou}{\mbox{\Ou}}

\newcommand{\p}{\:\:.}  
\newcommand{\PGL}{\mathrm{PGL}}
\newcommand{\PP}{\mathbbm{P}}

\newcommand{\Pic}{\mathrm{Pic}} 
\newcommand{\pic}{\mathrm{pic}}
\newcommand{\plus}{\oplus} \newcommand{\Plus}{\bigoplus}

\newcommand{\QQ}{\mathbbm{Q}}
 
\newcommand{\qq}{\forall\:}

\newcommand{\res}[1]{{\left | {}_{#1} \right.}}

\newcommand{\resp}{\emph{resp.}}

\newcommand{\SL}{\mathrm{SL}}

\newcommand{\Sp}{\mathrm{Sp}}
 
\newcommand{\Spin}{\mathrm{Spin}}
 \newcommand{\sta}{\stackrel}
\newcommand{\sub}{\subseteq}

 \newcommand{\tens}{\otimes} 

\newcommand{\tilda}{\widetilde} 

\newcommand{\tq}{\: : \:}

\newcommand{\vide}{\emptyset}

\newcommand{\ZZ}{\mathbbm Z}
\newcommand{\zzz}{{\mathbf{z}}}

\begin{document}

\title{Normalité projective des varétés magnifiques de rang 1\\
(Projective normality of rank one wonderful varieties)}
 \author{Alexis TCHOUDJEM\\
  Institut Camille Jordan\\
Universit\'e Claude Bernard Lyon I\\
Boulevard du Onze Novembre 1918\\
69622 Villeurbanne\\
FRANCE\\
tchoudjem@math.univ-lyon1.fr
} \date{Villeurbanne, le \today}

\maketitle

{\bf Résumé :} Soient $\call{L}$ et $\call{L'}$ deux faisceaux inversibles sur une variété projective $X$. On suppose que $\call{L}$ et $\call{L}'$ sont engendrés par leurs espaces de sections globales $\Gamma(\call{L})$ et $\Gamma(\call{L}')$. On démontre dans cet article que  le morphisme :
\[\Gamma(\call{L}) \tens \Gamma(\call{L}') \to \Gamma(\call{L} \tens \call{L}')\]
est surjectif, dans le cas où $X$ est une variété magnifique de rang $1$. En particulier,  le c\^one au-dessus d'une variété magnifique $X$ de rang $1$ défini par un  faisceau inversible très ample est toujours normal.

\vskip 1cm 

{\bf Abstract :} Let $\call{L}$ and $\call{L}'$ be two invertible sheaves over a projective variety $X$. We suppose that $\call{L}$ and $\call{L}'$ are generated by their global section spaces $\Gamma(\call{L})$ and $\Gamma(\call{L}')$. We prove in this article that the morphism :

\[\Gamma(\call{L}) \tens \Gamma(\call{L}') \to \Gamma(\call{L} \tens \call{L}')\]
is surjective, in the case where $X$ is a rank one wonderful variety. In particular, the cone over a rank one wonderful variety defined by a very ample invertible sheaf is always normal.

\section*{Introduction}

D'après \cite{chima}, si $X$ est une variété projective obtenue comme la compactification magnifique d'un espace symétrique $G/H$ (où $G$ est un groupe semi-simple adjoint et $H$ le sous-groupe des points fixes d'une involution $\theta$ de $G$), alors pour tous faisceaux inversibles et engendrés par leurs sections globales, le morphisme :

\begin{equation}\label{eq:un}
R \;: \; \Gamma(\call{L}) \tens \Gamma(\call{L}') \to \Gamma( \call{L} \tens \call{L}')
\end{equation}
induit par la restriction de $X \croi X$ à sa diagonale $\Delta_X$ est surjectif (on note $\Gamma(-)$ ou $\Gamma(X,-)$ l'espace des sections globales d'un faisceau), \cf aussi  \cite{kan}, pour les compactifications magnifiques des groupes semi-simples adjoints.

Les variétés magnifiques (définies par exemple dans \cite{lu}) sont une généralisation des compactifications magnifiques des espaces symétriques. Les variétés magnifiques les plus simples, après celles de rang $0$ qui sont les variétés de drapeaux (généralisées), sont les variétés magnifiques de rang $1$. Nous démontrons ici (\cf le théorème \ref{thm:prin}) que le morphisme $R$ de (\ref{eq:un}) est encore surjectif si $X$ est une variété magnifique de rang $1$.

La normalité du c\^one au-dessus de $X$ défini par un faisceau inversible très ample sur $X$ s'ensuit par un argument standard (\cf \cite[II Ex.5.14 d)]{ha}).

Avant de démontrer le résultat principal, nous rappelons dans la section \ref{sec:indu}, la définition des variétés magnifiques de rang $1$ et la notion d'induction parabolique. Pour une variété magnifique de rang $1$ donnée, nous aurons besoin d'une description des faisceaux inversibles sur $X$ et de l'espace de leurs sections globales (\cf la section \ref{sec:pic}) ainsi que de la grosse cellule de $X$, qui est un ouvert particulier isomorphe à un espace affine (\cf la section \ref{sec:gros}).

Dans toute la suite $\kk$ est un corps algébriquement clos de caractéristique nulle et $G$ est un groupe linéaire algébrique, semi-simple et simplement connexe sur $\kk$, d'algèbre de Lie $\goth g$. On choisit un sous-groupe de Borel $B$ de $G$ et $T\sub B$ un tore maximal. Soient $B^-$ le sous-groupe de Borel de $G$ opposé à $B$ relativement à $T$. On note $\Phi$ le système de racines de $(G,T)$, $\Phi^+$ l'ensemble des racines positives, $\Delta$ la base de $\Phi$ associée à $B$ et $W$ le groupe de Weyl de $(G,T)$. On notera $w_0$ l'élément de plus grande longueur de $W$.

\section{Induction parabolique}\label{sec:indu}

Soit $X$ une $G-$variété. Soit $P$ un sous-groupe parabolique de $G$.  S'il existe un morphisme $G-$équivariant $\pi : X \to G/P$, une $P-$variété $ Y$ et une immersion fermée $Y \to X$, $P-$équivariante d'image $\pi\inv(P/P)$, on dit que $X$ s'obtient par induction parabolique de $Y$. On le note :\[X= G\croi^PY\p\]

Une $G-$variété magnifique de rang $1$ est une $G-$variété projective lisse connexe avec deux $G-$orbites : une $G-$orbite ouverte : $X_G^0$ et une $G-$orbite fermée $F:=X\moins X_G^0$ de codimension $1$.

Par exemple, la variété $\PP^1 \croi \PP^1$ munie de l'action diagonale du groupe $G = \SL_2$ est une variété magnifique de rang $1$ (la diagonale et son complémentaires forment les deux $G-$orbites).

Une liste d'autres exemples est donnée par  la table $1$ de \cite{was}. On dira que ces variétés magnifiques de la table $1$ de \cite{was} et la $\SL_2-$variété $\PP^1 \croi \PP^1$ sont les variétés magnifiques {\it irréductibles} de rang $1$ (\cf aussi \cite{akh} et \cite{bri}).

Toutes les autres variétés magnifiques de rang $1$ s'obtiennent à partir des variétés irréductibles par induction parabolique :
 
\begin{lem}[{\cite[\S 2]{was} ou \cite{akh} ou \cite{bri}}]
Soit $X$ une $G-$variété magnifique de rang $1$. Alors, il existe un sous-groupe parabolique $P$ de $G$, contenant $B^-$, tel que :
\[X=G\croi^P\tilda{X}\]
où $\tilda{X}$ est une $P/P^r-$variété magnifique  irréductible.
\end{lem}

\section{Groupe de Picard}\label{sec:pic}
On rappelle quelques faits sur le groupe de Picard de $X$. Soit $H$ un sous-groupe d'isotropie d'un point de l'orbite ouverte de $X$. Quitte à changer $H$ en un de ses conjugués, on suppose que $BH/H$ est ouvert dans $G/H$. Soit $\lambda$ un poids dominant. Soit $V_\lambda$ le $G-$module simple de plus au poids $\lambda$. On suppose que $V_\lambda^{(H)}\neq 0$ \ie il existe $0\neq v_H \in V_\lambda$ un $H-$vecteur propre. Comme la variété $X$ est magnifique, le morphisme $gH \donne [g.v_H] \in \PP(V_\lambda)$ se prolonge en un morphisme :\[i_\lambda : X \to \PP(V_\lambda) \p\] On pose :\[\call{L}_\lambda :=i_\lambda^*(\call{O}(1))\]
c'est un faisceau inversible et $G-$linéarisé sur $X$ (qui {\it a priori} dépend non seulement de $\lambda$ mais aussi du choix d'un $H-$vecteur propre dans $V_\lambda$). Soit $l_\lambda \in V_\lambda^*$ un $B-$vecteur propre (de poids $-w_0\lambda$). On notera $\sigma_\lambda$ l'image de $l_\lambda \in V_\lambda^* = \Gamma(\PP(V_\lambda),\call{O}(1))$ dans $\Gamma(X,\call{L}_\lambda)$.

Soit maintenant  $D$ un diviseur premier et $B-$stable de $G/H$. Soit $\adh{D}$ son adhérence dans $X$. Soit $p : G \to G/H$, $ g \donne gH$. Il existe $f_D \in \kk[G]$ tel que $p\inv(D) = (f_D)$ comme diviseurs de $G$. Comme $D$ est $B-$stable, il existe un $G-$module simple $V_{\lambda_D}$, de plus haut poids $\lambda_D$ dominant, $l_{D} \in V_{\lambda_D}^*$ un $B-$vecteur propre et $v_H\in V_{\lambda_D}$ un $H-$vecteur propre tels que :
\[f_D = l_{D} \tens v_H\]
\ie : $f_D(g) = l_D(gv_H)$ pour tout $g \in G$. Alors, on peut montrer que :
\[\call{L}_{\lambda_D} \iso \call{O}_X(\adh{D}) \]
(en particulier, le faisceau $\call{O}_X(\adh{D})$ est engendré par ses sections globales).

De plus, le groupe de Picard de $X$, $\Pic (X)$, est un réseau dont les classes des faisceaux $\call{L}_{\lambda_D}$ forment une base. 
\[\]
On note $\zzz$ l'unique point fixe de $B^-$ dans $X$ et $Q$ son stabilisateur dans $G$. Alors $Q$ est un sous-groupe parabolique de $G$ contenant $B^-$. On note $P:=w_0Qw_0$ le sous-groupe parabolique de $G$ opposé à $Q$. Tout caract\`ere $\lambda$ de $P$ induit un faisceau inversible (et $G-$linéarisé) sur $G/Q$ : $\call{L}_{G/Q}(\lambda)$ (\cf \cite[I \S 5.8]{jan}).  On prendra pour convention que la fibre au-dessus de $Q/Q$,   $\call{L}_{G/Q}(\lambda)\res{Q/Q}$, est la droite affine  où le groupe $Q$ agit via le caractère $\lambda^*:=-w_0\lambda$.

Avec cette convention et en identifiant la $G-$orbite fermée $F$ de $X$ avec $G/Q$, si $\lambda$ est un poids dominant tel que $V_\lambda^{(H)}\neq 0$, alors $\call{L}_\lambda\res{F} = \call{L}_{G/Q}(\lambda)$.

\begin{lem}\label{lem:res}
Soit $F$ l'unique $G-$orbite fermée de $X$. On suppose que  $X$ n'est pas une induction parabolique de la variété $\PP^1\croi \PP^1$  munie de l'action diagonale de $\SL_2$. Alors  le morphisme :
\[\Pic (X) \to \Pic (F) \,,\; \call{L} \donne \call{L}\res{F}\]
est injectif de plus, un faisceau inversible $\call{L}$ sur $X$ est engendré par ses sections globales si et seulement si le faisceau $\call{L}\res{F}$ est engendré par ses sections globales sur $F$.
\end{lem}

\begin{dem}
On identifie la $G-$orbite fermée $F$ avec $G/Q$.

Soient $D_1,...D_p$ les diviseurs premiers $B-$stables de la $G-$orbite ouverte de $X$. Soient $\lambda_1,...,\lambda_p$ les caractères correspondants. Les classes des faisceaux $\call{L}_{\lambda_1},..., \call{L}_{\lambda_p}$ forment une base de $\Pic (X)$, et  pour tout $i$, $\call{L}_{\lambda_i} \res{G/Q} = \call{L}_{G/Q}(\lambda_i)$.

Si $X$ est une variété magnifique de rang $1$ irréductible, alors $p=1$ (et le lemme est vrai dans ce cas) ou, d'après le tableau $1$ de \cite[table 1]{was}, $p=2$ et $\lambda_1,\lambda_2$ sont des poids fondamentaux distincts et le lemme est encore vrai dans ce cas.  

Pour le cas général, on a $X = G\croi^{Q_1} X_1$ pour un certain sous-groupe parabolique $Q_1$ de $G$ contenant $B^-$ et $X_1$ une  $Q_1/Q_1^r-$variété magnifique irréductible de rang $1$.

Si on note $\pi : X \to G/Q_1$ et $i : X_1 \to X$, alors on a une suite exacte :
\[0 \to \Pic(G/Q_1) \sta{\pi^*}{\to} \Pic(X )\sta{i^*}{\to} \Pic (X_1) \to 0\p\]
En effet, si $\call{L}$ est un faisceau inversible sur $X$ tel que $\call{L}\res{X_1} \iso \call{O}_{X_1}$, alors on peut montrer que le faisceau $\pi_*\call{L}$ est inversible sur $G/Q_1$ et que $\pi^*\pi_*\call{L} \iso \call{L}$ ; pour la surjectivité de $i^*$, on remarque que la $B-$orbite ouverte $BQ_1/Q_1$ de $G/Q_1$ est isomorphe à un espace affine, que l'ouvert de $X$, $\Omega := \pi\inv(BQ_1/Q_1)$ est isomorphe à $BQ_1/Q_1 \croi X_1$ d'où un morphisme surjectif et un isomorphisme : $\Pic(X) \to \Pic(\Omega) \iso \Pic(X_1)$.

Donc si $\call{L}$ est un faisceau inversible sur $X$ tel que $\call{L}\res{G/Q} \iso \call{O}_{G/Q}$, alors la restriction de $\call{L}$ à $F_1$, l'orbite fermée de $X_1$ est triviale : $\call{L}\res{F_1} \iso \call{O}_{F_1}$. Donc $\call{L}\res{X_1} \iso \call{O}_{X_1}$ d'après le premier cas. Mais alors, $\call{L}$ est de la forme $\pi^*\call{M}$ pour un certain faisceau inversible $\call{M}$ sur $G/Q_1$. Si on note $\pi_1$ la restriction  de $\pi$ à $G/Q$, $\pi_1^ *\call{M} \iso \call{O}_{G/Q}$. 

\xymatrix{X_1\ar^i[r] & X \ar^\pi[r]& G/Q_1\\
F_1 \ar[r]\incl[u]& G/Q\incl[u]\ar^{\pi_1}[ur]& 
}

Or, $Q \sub Q_1$ et  $\pi_1$ est aussi la surjection canonique  $G/Q \to G/Q_1$. Comme le morphisme :
\[\Pic (G/Q_1) \sta{\pi_1^*}{\to} \Pic(G/Q)\]
est injectif, le faisceau inversible $\call{M}$ est trivial et $\call{L}= \pi^* \call{M}$ aussi.

Soit $D$ un diviseur de $X$ tel que le faisceau :
\[\call{O}_X(D)\res{F}\]
est engendré par ses sections globales. Le diviseur $D$ est équivalent à un diviseur $n_1\adh{D_1}+...+n_r\adh{D_r}$ pour certains entiers $n_i$. Pour montrer que $\call{O}_X(D)$ est engendré par ses sections globales, il suffit de montrer que $n_i \ge 0$ pour tout $i$.

Quitte à renuméroter les $D_i$, on peut supposer que pour un certain $r \le p$, $D_1,...,D_r \sub X \moins \Omega$ et que $D_i \cap \Omega \neq \vide$ si $i >r$. Soient $E_1,...,E_r$ les diviseurs premiers $B-$invariants de $G/Q_1$ tels que $\pi\inv (E_i) = \adh{D_i}$. Soient $\delta_i$, $i >r$ les diviseurs premiers $B \cap Q_1-$invariants de $X_1$ tels que $\adh{D_i}= \adh{BQ_1 /Q_1 \croi \delta_i}$.

Avec ces notations, les poids $\lambda_1^*,...,\lambda_r^*$ sont des poids fondamentaux (deux à deux distincts) qui se prolongent en des caractères de $Q_1$ et les poids $\lambda_{r+1}^*,...,\lambda_p^*$ sont des caractères triviaux sur $Q_1^r$.

Comme $\call{O}_X(D)\res{F}$ est engendré par ses sections globales, le poids $n_1\lambda_1^*+...+n_p\lambda_p^*$ est dominant. Donc :
\[\qq 1 \le i \le r,\; n_i \ge 0 \p\]

Notons $F_1$ la $Q_1-$orbite fermée de $X_1$. Le faisceau 
\[\call{O}_X(D)\res{F_1} \iso  \call{O}_{X_1}(\sum_{i>r} n_i {\delta_i})\res{F_1}\]
est engendré par ses sections globales donc, d'après le cas irréductible, le faisceau $\ma \call{O}_{X_1}(\sum_{i>r} n_i {\delta_i})$ est aussi engendré par ses sections globales d'où :
\[\qq i >r,\; n_i \ge 0\p\]
\end{dem}

{\it Remarque : si $X = \PP^1 \croi \PP^1$, $F=\Delta_{\PP^1} \iso \PP^1$ et le morphisme $\Pic X = \ZZ^2  \to \Pic F =\ZZ$ ne peut pas \^etre injectif de plus, dans ce cas, le faisceau $\call{O}(1) \cten \call{O}(-1)$ n'est pas engendré par ses sections globales alors que sa restriction $\call{O}(1) \cten \call{O}(-1) \res{\Delta_{\PP^1}} \iso \call{O}_{\PP^1}$ l'est. 

Excepté ce cas et ses inductions paraboliques, on obtient en particulier que pour tout poids dominant $\lambda$, le faisceau inversible $\call{L}_\lambda$ défini précedemment ne dépend que de $\lambda$ car $\call{L}_\lambda \res{F} \iso \call{L}_{G/Q}(\lambda)$.} 

\vskip 2em

Si $X$ n'est pas une induction parabolique de $\PP^1 \croi \PP^1$, on en déduit :

\begin{cor}

Si $\call{L}$ est un faisceau inversible sur $X$ engendré par ses sections globales, alors il existe $\lambda$ un poids dominant tel que \[(*)\;\; V_\lambda^{(H)} \neq 0\] et $\call{L} \iso \call{L}_\lambda$. Si $\lambda,\mu$ sont deux poids dominants qui vérifient $(*)$, alors $\lambda+\mu$ vérifie aussi $(*)$ et $\call{L}_\lambda \tens \call{L}_\mu \iso \call{L}_{\lambda + \mu}$.
\end{cor}

\begin{dem}
Soient $\lambda,\mu$ des poids dominants qui vérifient $(*)$. Soient $v_H \in V_\lambda$ et $v'_H \in V_\mu$ deux $H-$vecteurs propres. Comme $BH/H$ est un ouvert dense de $G/H$, le vecteur $v_H$ a une composante non nulle de $T-$poids $w_0\lambda$ et le vecteur $v'_H$ une composante non nulle de $T-$poids $w_0\mu$. Donc le $H-$vecteur propre $v_H \tens v'_H$ a une composante non nulle de $T-$poids $w_0(\lambda + \mu)$ dans $V_\lambda \tens V_\mu$. Donc en considérant la projection $G-$équivariante : $V\lambda \tens V_\mu \to V_{\lambda+\mu}$, on trouve un $H-$vecteur propre dans $V_{\lambda+\mu}$.

Or :\[\left(\call{L}_\lambda\tens \call{L}_\mu \right)\res{G/Q} \iso  \call{L}_\lambda\res{G/Q} \tens \call{L}_\mu \res{G/Q}\]
\[\iso \call{L}_{G/Q}(\lambda) \tens \call{L}_{G/Q}(\mu)  \]
\[\iso \call{L}_{G/Q}(\lambda + \mu) \]
\[\iso \call{L}_{\lambda+\mu}\res{G/Q}\p \]
Donc $\call{L}_{\lambda+\mu} \iso \call{L}_\lambda\tens \call{L}_\mu$.

Soit $\call{L}$ un faisceau inversible engendré par ses sections globales.

Notons $D_1,...,D_p$ les diviseurs premiers $B-$invariants de $G/H$. Soient $\lambda_1,...,\lambda_p$ les poids correspondants (tels que $\call{L}_{\lambda_i} \iso \call{O}_X(\adh{D_i})$ pour tout $i$). Les poids $\lambda_i$ sont dominants.. De plus, il existe des entiers $n_1,...,n_p$ tels que :\[\call{L} \iso \call{L}_{\lambda_1}^{n_1}\tens ... \tens \call{L}_{\lambda_p}^{n_p}\]
\[\iso  \call{L}_{n_1 \lambda_1 +...+n_p\lambda_p}\]
  d'après la première partie de cette démonstration. Or, $\call{L}\res{G/Q}\iso \call{L}_{G/Q}(n_1 \lambda_1 +...+n_p\lambda_p)$ est aussi engendré par ses sections globales donc  $n_1 \lambda_1 +...+n_p\lambda_p $ est dominant.
\end{dem}

{\bf Notations :} On notera $\pic^+(X)$ les caractères $\lambda$ dominants tels que $V_\lambda^{(H)}\neq 0$ et pour tout $\lambda \in \pic^+(X)$, $\call{L}_\lambda$ un faisceau inversible et $G-$linéarisé sur $X$ correspondant (\ie défini comme au début de la section \ref{sec:pic}) (si $X$ est une induction parabolique de $\PP^1 \croi \PP^1$ et seulement dans ces cas, il peut y avoir deux faiscaux $\call{L}_\lambda$ non isomorphes).

\section{Grosses cellules}\label{sec:gros}

On note $P^u$ le radical unipotent du sous-groupe parabolique $P$ de $G$ et $L$ son sous-groupe de Levi contenant le tore $T$ (on rappelle que, relativement au tore $T$, $P$ est le sous-groupe parabolique opposé à $Q$ le stabilisateur de $\zzz$ l'unique point fixe de $B^-$ dans $X$). 

\begin{pro}[{\cite[th. 1.4 et cor. 1.5]{blv}}]
Il existe un voisinage ouvert affine $P-$stable de $\zzz$ dans $X$, noté $X_0$, une sous-variété fermée $A$ de $X_0$, $L-$ stable,  $L-$isomorphe à une droite affine (qui est une représentation de $L$) et un isomorphisme de $P-$variétés algébriques :
\[P^u \croi A \to X_0 \;,\; (u,a) \donne u.a\]
où $P$ agit sur le membre de gauche par :
\[(xl).(u,a) := (xlul\inv,l.a)\]
pour tous $x,u \in P^u$, $l\in L$, $a \in A$.
\end{pro}

L'ouvert de $X_0$ est la grosse cellule de $X$ (on peut l'obtenir comme la cellule ouverte d'une décomposition cellulaire de Bialynicki-Birula de $X$ associée à un sous-groupe à un paramètre dominant générique de $T$). 

{\it Remarque :} Comme $X$ n'a qu'une seule orbite fermée : $F = G.\zzz$,  et comme $\displaystyle X \moins \bigcup_{g \in G }gX_0$ est un fermé de $X$ ne contenant pas $\zzz$, on a forcément $\ma \bigcup_{g\in G} gX_0 = X$.

\`A partir de maintenant et jusqu'à la fin de cette section \ref{sec:gros}, on suppose que $X$ n'est pas une induction parabolique de $\PP^1 \croi \PP^1$.

Soient $\lambda$ un poids dominant tel que $V_\lambda^{(H)} \neq 0$ , $\call{L}_\lambda$ le faisceau inversible sur $X$ correspondant et $\sigma_\lambda \in \Gamma(X,\call{L}_\lambda)$ sa section de plus haut poids (définie à multiplication par un scalaire près). Alors on a :
\[\Gamma (X_0,\call{L}_\lambda) = \kk[X_0]\sigma_\lambda \p\]

De plus, comme d'après la remarque ci-dessus, $X$ est l'unique ouvert $G-$stable de $X$ contenant $X_0$, l'espace des sections globales  $\Gamma(X,\call{L}_\lambda)$ est le plus grand sous-$G-$module rationnel du $\goth g-$module $\Gamma(X_0,\call{L}_\lambda)$.

On a : $\kk[X_0] = \kk[P^u]\tens_\kk\kk[A]$. Soit $\y$ un générateur de l'algèbre $\kk[A]$ qui est un $L-$vecteur propre. Le poids de $y$ est $-\gamma$ où $\gamma$ est la {\it racine sphérique} de $X$ et on a :
\[\kk[X_0] = \plus_{m \ge 0} \kk[P^u]\y^m\]
et les fonctions $\y^m$ sont, à multiplication par une constante non nulle près, les $B-$vecteurs propres dans le $B-$module  $\kk[X_0]$. 

 Le lemme suivant est une adaptation d'un résultat plus général sur les variétés sphériques quelconques, \cf \cite[pro. du \S 3.3]{bripic}. On en donne néanmoins une démonstration {\it ad hoc}. On suppose que $X$ n'est pas une induction parabolique de $\PP^1 \croi \PP^1$. On rappelle que $\lambda^*=-w_0 \lambda$.

\begin{lem}\label{lem:dec}
Soit $\lambda$ dominant tel que $V_\lambda^{(H)}\neq 0$. Si $m \ge 0$ est un entier tel que le poids $\lambda^* -m\gamma$ est dominant, alors le $\goth g-$sous-module de $\Gamma(X_0,\call{L}_\lambda)$ engendré par $y^m\sigma_\lambda$, que nous noterons 
\[U(\goth g)\y^m\sigma_\lambda\;,\]
 est un sous-$G-$module irréductible de l'espace des sections globales $\Gamma(X, \call{L}_\lambda)$ et on a la décomposition suivante en somme directe de $G-$modules irréductibles :
\[\Gamma(X,\call{L}_\lambda) = \Plus_{m \ge 0 \atop \lambda^* - m \gamma \mathrm{\; dominant}}U(\goth g) y^m\sigma_\lambda\p\]

\end{lem}

\begin{dem}
Le $G-$module $\Gamma(X,\call{L}_\lambda)$ est engendré par ses $B-$vecteurs propres. Or les $B-$vecteurs propres de l'espace  $\Gamma(X_0,\call{L}_\lambda)$ sont de la forme $y^m\sigma_\lambda$ (à un scalaire près). Or dans un $G-$module les $B-$vecteurs propres sont de poids dominants donc :\[
\Gamma(X,\call{L}_\lambda) \sub \Plus_{m \ge 0 \atop \lambda^* - m \gamma \mathrm{\; dominant}}U(\goth g) y^m\sigma_\lambda\p\]

Pour l'autre inclusion, on remarque que le faisceau d'idéaux $\call{I}_F$ qui définit l'orbite fermée est un faisceau inversible $G-$linéarisé et que le tore $T$ agit sur la fibre $\call{I}_F\res{\zzz}$ via le caractère $-\gamma$. On en déduit donc que pour tout $m \ge 0$ tel que $\lambda^*-m\gamma$ est dominant, le faisceau $\call{L}_\lambda \tens \call{I}_F^m$ est engendré par ses sections globales (\cf le lemme \ref{lem:res}) et que  :
\[\call{L}_\lambda \tens \call{I}_F^m \iso \call{L}_{\lambda -m \gamma^*}\p\]

D'où un morphisme injectif de faisceaux $G-$linéarisés :
\[\call{L}_{\lambda - m \gamma^* } \to \call{L}_\lambda\p\]
Mais alors, le $G-$module $\Gamma ( X , \call{L}_{\lambda  - m \gamma^* })$ a une image non nulle dans $\Gamma ( X , \call{L}_{\lambda  })$. Or, l'image de $\sigma_{\lambda-m\gamma^*}$ est $y^m\sigma_\lambda$ (à multiplication par un scalaire non nul près)  donc $y^m\sigma_\lambda \in \Gamma ( X , \call{L}_{\lambda  })$.
\end{dem}

\section{Surjectivité de la multiplication}

\begin{thm}\label{thm:prin}
Soient $\call{L},\call{L'}$ deux faisceaux inversibles engendrés par leurs sections globales sur une $G-$variété magnifique de rang $1$. Alors le morphisme :
\[\Gamma(X,\call{L}) \tens \Gamma(X,\call{L'}) \to \Gamma(X,\call{L}\tens \call{L'})\]
est surjectif.
\end{thm}

\begin{dem}

On notera $\gamma $ la racine sphérique de $X$.

{\it Supposons pour commencer que $X$ n'est pas une induction parabolique de $\PP^1 \croi \PP^1$}.

Soient $\lambda,\mu$ des poids dominants tels que :\[ V_\lambda^{(H)}\neq 0\,,\; V_\mu^{(H)} \neq 0,\; \call{L}\iso \call{L}_\lambda \et \call{L}' \iso \call{L}_\mu \p\]

On a alors $\call{L}\tens \call{L}' \iso \call{L}_{\lambda +\mu}$.
Notons $R_{\lambda,\mu}$ ou $R$ s'il n'y a pas d'ambigu\"{i}té le morphisme :
 \[\Gamma ( X , \call{L}_{\lambda  }) \tens \Gamma ( X , \call{L}_{\mu  })\to \Gamma ( X , \call{L}_{\lambda  +\mu}) \p\]

{\bf Cas où $\mathbf X$ est une variété magnifique irréductible}

Dans ce cas, d'après \cite{akh}, la variété $X$ est soit isomorphe à une variété de drapeaux  (un espace projectif $\PP^n$, une quadrique $Q(n)$ d'équation $x_0^2=x_1^2+...+x_n^2$ dans un espace projectif $\PP^n$,  une variété grassmannienne $\mathrm{Gr}(2n,2)$ (variété des $2-$plans dans $\kk^{2n}$), un produit d'espace projectis $\PP^n \croi \PP^n$, ou une variété homogène pour le groupe adjoint de type $E_6$) soit isomorphe à une des variétés numérotées $9B$, $9C$, $15$ dans le tableau $1$ de \cite[tab.1]{was}.

--- Lorsque $X$ est une variété de drapeaux, le morphisme :
\[\Gamma(X,\call{L}) \tens \Gamma(X,\call{L'}) \to \Gamma(X,\call{L}\tens \call{L'})\] est surjectif  car c'est un morphisme non nul et l'espace $\Gamma(X,\call{L}\tens \call{L'})$ est une représentation irréductible de l'algèbre de Lie du groupe des automorphismes de $X$ (théorème de Borel-Weil).

--- Si $X$ est la variété numéro $9B$ (\resp\ $9C$)  dans le tableau $1$ de \cite[tab.1]{was}, alors il existe $n \ge 2$ tel que $G = \Spin_{2n+1}$ (\resp\ $G = \Sp_{2n}$). Notons $\alpha_1,...,\alpha_n$ les racines simples et $\omega_1,...,\omega_n$ les poids fondamentaux correspondants. Il existe $p,q \ge 0$ tels que $\lambda = p\omega_1,\mu=q\omega_1$, (\resp\ $\lambda = p\omega_2,\mu=q\omega_2$) et on a : $\gamma = \omega_1$ (\resp\ $\gamma = \omega_2$). Remarquons que dans les deux cas, on a $\lambda^* = \lambda$ et $\mu^*=\mu$.

Soit $m \ge 0$ tel que $\lambda + \mu -m \gamma$ est dominant. Cette condition est équivalente à :\[0\le m \le p+q\p \]

On pose $m_1 := \max\{m-q,0\} $ et $m_2:=m-m_1$. On a alors : $0\le m_1 \le p$, $0\le m_2 \le q$ et $m = m_1+m_2$ \cad : $\lambda - m_1 \gamma$, $\mu -m_2\gamma$ dominants. Donc d'après le lemme \ref{lem:dec} :
\[y^{m_1}\sigma_\lambda \in \Gamma(X,\call{L}_\lambda) \et y^{m_2}\sigma_\mu \in \Gamma(X,\call{L}_\mu) \p\]
Or, $R(y^{m_1}\sigma_\lambda \tens y^{m_2}\sigma_\mu) = y^{m}\sigma_{\lambda+\mu}$.

Donc l'image de $R$ est un sous $G-$module de $\Gamma(X,\call{L}_{\lambda+\mu})$ qui contient  tous les $y^{m}\sigma_{\lambda+\mu}$ tels que $m\ge 0$ et $\lambda +\mu - m \gamma$ est dominant. On en déduit donc que $R$ est surjectif.

--- Si $X$ est la variété numéro $15$ dans le tableau $1$ de \cite[tab.1]{was}, alors $G$ est le groupe semi-simple et simplement connexe de type $G_2$. On prend les notations standards :  $\alpha_1, \alpha_2$ pour les racines simples, $\omega_1,\omega_2$ pour les poids fondamentaux correspondants. On a $\gamma = \omega_2 - \omega_1$ et il existe des entiers $p,q,p',q'\ge 0$ tels que $\lambda = p\omega_1 + q\omega_2, \mu = p'\omega_1+q'\omega_2$. On a : $\lambda^*=\lambda$, $\mu^* =\mu$.

Soit $m\ge 0$ un entier tel que $\lambda+\mu - m \gamma $ est dominant.

Comme  $\lambda+\mu - m \gamma =(p+p'+m)\omega_1 +(q+q'-m)\omega_2$, on a : $0 \le m \le q+q'$. 

On pose : $m_1:=\max\{m-q',0\}$ et $m_2:= m-m_1$. On a $0\le m_1 \le q$, $0\le m_2 \le q'$ et $m= m_1+m_2$. Donc : $\lambda - m_1 \gamma$ et $\mu -m_2\gamma$ sont dominants,
\[y^{m_1}\sigma_{\lambda} \et y^{m_2}\sigma_\mu \in \Gamma(X,\call{L}_\lambda) \p\]
Or : $R(y^{m_1}\sigma_{\lambda} \tens y^{m_2}\sigma_\mu) = y^{m_1+m_2}\sigma_{\lambda+\mu}$. 

On en déduit donc que $R$ est encore surjectif dans ce cas.

{\bf Cas général}

Comme on l'a rappelé au début, il existe un sous-groupe parabolique $Q$ contenant $B^-$ tel que $X = G\croi^{Q}X_1$ pour une certaine $Q/Q^r-$variété magnifique irréductible. Soit $p : X \to G/Q$ la projection associée.

Notons $L$ le sous-groupe de Levi de $G$ contenant $T$ tel que $Q =Q^uL$. 
On notera $\Delta_1$ la base du système de racines associé à $(L,T)$ contenue dans $\Delta$.

Soit $B\cap Q /Q^r \sub P_1 \sub Q/Q^r$ le sous-groupe parabolique associé à la variété magnifique $X_1$. La grosse cellule $(X_1)_0$ de $X_1$ vérifie :
\[(X_1)_0 \iso P_1^u \croi A_1\]
pour une certaine sous-variété fermée $P_1/P_1^u-$stable $A_1$ de $X_1$. Si on note $Q^+$ le sous-groupe parabolique de $G$ opposé à $Q$ (relativement à $T$), alors la grosses cellule $X_0$ de $X$ est isomorphe à
\[(Q^+)^u \croi (X_1)_0 \iso (Q^+)^u \croi  P_1^u \croi A_1 \p\]

Soit $y \in \kk[X_0]$ un $P_1/P_1^u-$vecteur propre tel que  $\kk[A_1]=\kk[y]$. Alors $y$ est de poids $\gamma$, la racine sphérique de $X$. En particulier, $\gamma$ est un caractère de $T/Q^r$ et donc comme $\ma Q^r=\left( \bigcap_{\delta \in \Delta_1}\ker \delta\right)^\circ$, $\gamma \in \QQ\Delta_1$.

On a une décomposition des poids $\lambda^*$ et $\mu^*$ :
\[\lambda^* =\underbrace{\sum_{\alpha \in \Delta_1}\lambda_\alpha \omega_\alpha}_{=:\lambda_1^*} + \underbrace{\sum_{\alpha \in \Delta\moins \Delta_1}\lambda_\alpha\omega_\alpha}_{=:\lambda_2^*}\]
\[\mu^* =\underbrace{\sum_{\alpha \in \Delta_1}\mu_\alpha \omega_\alpha}_{=:\mu_1^*} + \underbrace{\sum_{\alpha \in \Delta\moins \Delta_1}\mu_\alpha\omega_\alpha}_{=:\mu_2^*} \]
 pour certains poids dominants $\lambda_1,\lambda_2,\mu_1,\mu_2$.

Comme $\langle \lambda_2^*, \alpha^\ch\rangle = 0$ pour tout $\alpha \in \Delta_1$, le caractère $\lambda_2^*$ se prolonge en un caractère de $Q$. Donc $\lambda_2 \in \pic^+(X)$ et $\call{L}_{\lambda_2} \iso p^*\call{L}_{G/Q}(\lambda_2)$. On en déduit aussi que $\lambda_1\in \pic^+(X)$ avec $\call{L}_{\lambda_1} \iso \call{L}_\lambda \tens \call{L}_{\lambda_2}\inv$.   

De m\^eme, $\mu_2^*$ se prolonge en un caractère de $Q$ et $\mu_1 \in \pic^+(X)$.

On pose $L':=(L,L)$ et on choisit un tore maximal $T'$ de $L'$ contenu dans $T$.

Si $\nu$ est un caractère de $T$, on notera $\nu'$ sa restriction à $T'$.

Si  $\lambda_2=\mu_2=0$, alors, le morphisme de restriction :
\[\Gamma(X,\call{L}_{\lambda_1}) \to \Gamma(X_1,\call{L}_{\lambda_1}\res{X_1})\]
est surjectif.

En effet, si $\sigma_{\lambda_1}\in \Gamma(X,\call{L}_{\lambda_1})$ est un vecteur propre de plus haut poids $\lambda_1$ et si on pose :$\sigma'_{\lambda_1}:=\sigma_{\lambda_1}\res{X_1}$ et $y':=y\res{(X_1)_0}$, alors daprès le lemme \ref{lem:dec}, on a :
\[\Gamma(X,\call{L}_{\lambda_1}) = \Plus_{m \ge 0 \atop \lambda_1^* - m\gamma \mathrm{\, dominant}}U(\goth g)y^m\sigma_{\lambda_1} \]
\[\Gamma(X_1,\call{L}_{\lambda'_1}) = \Plus_{m \ge 0 \atop (\lambda^*_1)' - m\gamma' \mathrm{\, dominant}}U(\goth l'){y'}^m\sigma'_{\lambda_1} \]
où $\goth l'$ est l'algèbre de Lie de $L'$.

Or, comme $\gamma \in  \QQ\Delta_1$ et comme $\lambda_1^* \in \langle \omega_\delta \tq \delta \in \Delta_1\rangle$, $\lambda_1^* -m\gamma$ est dominant si et seulement si $(\lambda_1^*)' - m\gamma'$ est dominant. D'où la surjectivité.

De m\^eme, le morphisme de restriction :
\[\Gamma(X,\call{L}_{\mu_1}) \to \Gamma(X_1,\call{L}_{\mu_1}\res{X_1})\]
est surjectif.

Or, d'après le cas particulier des variétés magnifiques irréductibles, le morphisme :
\[\Gamma(X_1,\call{L}_{\lambda_1} )\tens\Gamma(X_1,\call{L}_{\mu_1}) \to \Gamma(X_1,\call{L}_{\lambda_1 + \mu_1})\]

est surjectif.

Pour montrer la surjectivité de $R$, il suffit de montrer que si $m \ge 0$ est tel que $\lambda_1^*+\mu_1^* - m \gamma$ est dominant, alors :\[y^m\sigma_{\lambda_1+\mu_1} \in \im R\p\]

Considérons donc  $m \ge 0$ tel que  $\lambda_1^*+\mu_1^* - m \gamma$ est dominant. On a alors $\lambda'_1+\mu'_1-m\gamma'$ dominant.

Puisque le morphisme $R$ commute avec la restriction à $X_1$ :

\[\xymatrix{
\Gamma(X,\call{L}_{\lambda_1}) \tens \Gamma(X,\call{L}_{\mu_1}) \ar@{->>}[d]\ar^-{R}[r] & \Gamma(X,\call{L}_{\lambda_1+\mu_1})\ar@{->>}[d]\\
\Gamma(X_1,\call{L}_{\lambda_1}) \tens \Gamma(X_1,\call{L}_{\mu_1}) \ar@{->>}[r] & \Gamma(X_1,\call{L}_{\lambda_1+\mu_1})
}\]

il existe $\sigma \in \Gamma(X,\call{L}_{\lambda_1}) \tens \Gamma(X,\call{L}_{\mu_1})$ tel que :
\[R(\sigma)\res{X_1} = y^m\sigma_{\lambda_1+\mu_1}\res{X_1} \p\]

Dans le diagramme ci-dessus, les morphismes sont $\diag(T)-$équivariants (car $X_1$ est une sous-variété $Q-$stable de $X$). On peut donc choisir $\sigma$ comme un $\diag(T)-$vecteur propre de poids $\lambda_1^*+\mu_1^* - m\gamma$.

Nous allons voir qu'alors  $R(\sigma) = y^m\sigma_{\lambda_1+\mu_1}$.

En effet, $R(\sigma) = h \sigma_{\lambda_1+\mu_1}$ pour un certain $T-$vecteur propre $h\in \call{O}(X_0)$ tel que $h\res{(X_1)_0} = y^m\res{(X_1)_0}$.

En particulier, $h =y^m +h'$ où $h'$ est un $T-$vecteur propre de $\call{O}(X_0)$ tel que $h'\res{(X_1)_0} = 0$.

Or on a des isomorphismes $T-$équivariants de variétés :
\[X_0 \iso (Q^+)^u \croi P_1^u \croi A_1 \et (X_1)_0 \iso \{1\} \croi P_1^u \croi A_1\]
donc si on note $u_1,...,u_p$ un système de coordonnées $T-$équivariantes de $(Q^+)^u$, $v_1,...,v_q$ un système de coordonnées $T-$équivariantes de $P_1^u$, alors $h'$ est une combinaison $\kk-$linéaire de mon\^omes :
\[u_1^{m_1}...u_p^{m_p}v_1^{n_1}...v_q^{n_q}y^k\]
où $m_1,...,m_p,n_1,...,n_q,k$ sont des entiers $\ge 0$ et au moins un des $m_i$ est non nul (car $h'\res{(X_1)_0} = 0$).

Mais si $-\alpha_1,...,-\alpha_p$, $-\beta_1,...,-\beta_q$ sont les $T-$poids des coordonnées $u_i,v_j$, alors :
\[\left\{ \alpha_1,....,\alpha_p \right\} = \Phi^+\moins \langle \Delta_1 \rangle \et \left\{\beta_1,...,\beta_q\right\} \sub \langle \Delta_1\rangle \p\]  

Comme $h'$ est de $T-$poids $-m\gamma$, les mon\^omes qui apparaissent dans la décomposition de $h'$ doivent vérifier :
\[- m_1\alpha_1-...-m_p\alpha_p -n_1\beta_1 -...-n_q\beta_q -k\gamma =  - m\gamma\]
\[\equi (m-k)\gamma -n_1\beta_1-...-n_q\beta_q = m_1\alpha_1+...+m_p\alpha_p\]

ce qui est impossible : le membre de gauche est dans $\QQ\Delta_1$ alors que  dans le membre de droite appara\^it au moins un coefficient non nul suivant une racine simple $\delta \in \Delta \moins \Delta_1$ car tous les $\alpha_i$ sont dans $\Phi^+\moins \Delta_1$ et $m_i>0$ pour au moins un indice $i$.

On en déduit que $h'=0$.

Pour le cas général si $\lambda_2$ ou $\mu_2 $ n'est pas forcément nul, on utilise le diagramme commutatif suivant :
\[\xymatrix{
\Gamma(\call{L}_{\lambda_1})\tens\Gamma(\call{L}_{\lambda_2})\tens \Gamma(\call{L}_{\mu_1})\tens\Gamma(\call{L}_{\mu_2}) \ar[r]^-*+{\scriptscriptstyle R_{\lambda_1,\mu_1}\tens R_{\lambda_2,\mu_2}\atop \; } \ar_-{R_{\lambda_1,\lambda_2}\tens R_{\mu_1,\mu_2}}[d] & \Gamma(\call{L}_{\lambda_1 + \mu_1})\tens\Gamma(\call{L}_{\lambda_2+\mu_2})\ar^-{R_{\lambda_1+\mu_1,\lambda_2+\mu_2}}[d]\\
\Gamma(\call{L}_{\lambda_1+\lambda_2}) \tens \Gamma(\call{L}_{\mu_1+\mu_2})\ar^-{R_{\lambda,\mu}}[r]& \Gamma(\call{L}_{\lambda_1+\lambda_2+\mu_1+\mu_2})
}
\p\]

Pour montrer que $R_{\lambda,\mu}$ est surjectif, il suffit donc de montrer que \[R_{\lambda_1,\lambda_2}, R_{\mu_1,\mu_2}, R_{\lambda_2,\mu_2}, R_{\lambda_1+\mu_1,\lambda_2+\mu_2}\] sont surjectifs (car on sait déjà que $R_{\lambda_1,\mu_1}$ est surjectif).

Or, pour tout $\lambda \in \pic^+(X)$ et tout caractère dominant $\mu_2$ de $Q$, le morphisme $R_{\lambda,\mu_2}$ est surjectif. En effet, si $m \ge0$ et si $\lambda^* +\mu_2 - m\gamma$ est dominant, alors $\lambda^*-m\gamma$ est aussi dominant  puisque l'on a  :

\[  \delta \in \Delta_1 \impliq \langle \mu_2 , \delta^\ch\rangle = 0 \et \delta \in \Delta \moins \Delta_1 \impliq  \langle \gamma , \delta^\ch\rangle \le 0 \p\]

Donc dans ce cas, $y^m\sigma_\lambda \in \Gamma(X,\call{L}_{\lambda})$ et $y^m\sigma_{\lambda+\mu_2} = R_{\lambda,\mu_2}(y^m\sigma_\lambda \tens \sigma_{\mu_2})$.   

\vskip 0.5cm

{\it On termine avec le cas particulier où $X$ est une induction parabolique de $\PP^1 \croi \PP^1$}.

Il existe $Q$ un sous-groupe parabolique de $G$ contenant $B^-$ tel que $Q/Q^r \iso \PGL_2$ ou $\SL_2$ et $\pi : X \to G/Q$ un morphisme $G-$équivariant tel que :
\[\pi\inv (Q/Q) \iso \PP^1 \croi \PP^1\]
(isomorphisme $\PGL_2-$équivariant (pour l'action diagonale de $\PGL_2$ sur $\PP^1 \croi \PP^1$)). On notera $X_1 := \pi\inv(Q/Q)$.

Notons $f$ et $m$ les morphismes suivants :
\[\xymatrix{
G \croi X_1 \ar^m[d]\ar^f[r]& G/Q \croi X_1\\
X 
}\] 
$f : (g,x) \donne (gQ,x)$, $m:(g,x) \donne g.x$. Pour tout $\lambda$ caractère de $Q$ et tous $k_1,k_2$ entiers, on note $\call{L}_{\lambda,k_1,k_2}$ le faisceau $m_*f^*\left(\call{L}_{G/Q}(\lambda)\cten \call{O}_{\PP^1 \croi \PP^1}(k_1,k_2)\right)$ sur $X$.

Les faisceaux $\call{L}_{\lambda,k_1,k_2}$ sont inversibles (il suffit de le vérifier au-dessus de l'ouvert $\pi\inv(BQ/Q)$ de $X$ qui est isomorphe à $BQ/Q \croi \PP^1 \croi \PP^1$) et tous les faisceaux inversibles sur $X$ sont de cette forme à isomorphisme près à cause de la suite exacte :

\[0 \to \Pic(G/Q) \sta{\pi^*}{\to} \Pic(X ) {\to} \Pic (X_1) \to 0\]
(\cf la démonstration du lemme \ref{lem:res}).

Soient $\call{L}$ et $\call{L'}$ des faisceaux inversibles sur $X$ engendrés par leurs sections globales. Soient $\lambda,\mu$ des caractères de $Q$, $k_1,k_2,l_1,l_2$ des entiers tels que :
\[\call{L} \iso \call{L}_{\lambda,k_1,k_2}\et \call{L'} \iso \call{L}_{\mu,l_1,l_2} \p\]

On a alors :
\[\call{L} \tens \call{L'} \iso \call{L}_{\lambda+\mu,k_1+k_2,l_1+l_2}\p\]

Puisque $\call{L}$ et $\call{L'}$ sont engendrés par leurs sections globales, les poids $\lambda,\mu$ sont des caractères de $T$ dominants (relativement à $B$) et on a:  $k_1,k_2,l_1,l_2 \ge 0$.

Si $s \ge 0$ est entier, on note $\kk[X,Y]_s$ est l'espace des polyn\^omes homogènes de degré $s$ en les deux variables $X,Y$. On a :
\[\Gamma(X,\call{L}) \iso V_{\lambda^*} \tens_{\kk} \kk[X,Y]_{k_1} \tens_{\kk}\kk[X,Y]_{k_2}\]
\[\Gamma(X,\call{L}) \iso V_{\mu^*} \tens_{\kk} \kk[X,Y]_{l_1} \tens_{\kk}\kk[X,Y]_{l_2}\]
\[\Gamma(X,\call{L} \tens \call{L'}) \iso V_{\lambda^*+\mu^*} \tens_{\kk} \kk[X,Y]_{k_1+l_1} \tens_{\kk}\kk[X,Y]_{k_2+l_2} \p\]

La surjectivité du morphisme :
\[\Gamma(X,\call{L})\tens \Gamma(X,\call{L'}) \to \Gamma(X,\call{L}\tens \call{L'})\] résulte alors directement de la surjectivité des morphismes :

\[V_{\lambda^*} \tens V_{\mu^*} \to V_{\lambda^* + \mu^*}\]
et \[\kk[X,Y]_{k_i} \tens \kk[X,Y]_{l_i} \to \kk[X,Y]_{k_i+l_i}\p\]
\end{dem}
\bibliographystyle{plain-fr}
\bibliography{biblio}

\begin{thebibliography}{10}
\expandafter\ifx\csname fonteauteurs\endcsname\relax
\def\fonteauteurs{\scshape}\fi

\bibitem{akh}
D.~\bgroup\fonteauteurs\bgroup Akhiezer\egroup\egroup{} :
\newblock Equivariant completions of homogeneous algebraic varieties by
  homogeneous divisors.
\newblock {\em Ann. Global Anal. Geom.}, 1(1)\string:\penalty500\relax 49--78,
  1983.

\bibitem{bripic}
M.~\bgroup\fonteauteurs\bgroup Brion\egroup\egroup{} :
\newblock Groupe de picard et nombres caractéristiques des variétés sphériques.
\newblock {\em Duke Math. J.}, 58(2)\string:\penalty500\relax 397--424, 1989.

\bibitem{bri}
M.~\bgroup\fonteauteurs\bgroup Brion\egroup\egroup{} :
\newblock On spherical varieties of rank one.
\newblock {\em in Group actions and invariant theory ({M}ontreal, PQ, 1988)},
  pages 31--41, 1989.

\bibitem{blv}
M.~\bgroup\fonteauteurs\bgroup Brion\egroup\egroup{},
  D.~\bgroup\fonteauteurs\bgroup Luna\egroup\egroup{} et Th.
  \bgroup\fonteauteurs\bgroup Vust\egroup\egroup{} :
\newblock Espaces homogènes sphériques.
\newblock {\em Invent. Math.}, 84(3)\string:\penalty500\relax 617--632, 1986.

\bibitem{chima}
R.~\bgroup\fonteauteurs\bgroup Chirivi\egroup\egroup{} et
  A.~\bgroup\fonteauteurs\bgroup Maffei\egroup\egroup{} :
\newblock Projective normality of complete symmetric varieties.
\newblock {\em Duke Math. J.}, 122(1)\string:\penalty500\relax 93--123, 2004.

\bibitem{ha}
R.~\bgroup\fonteauteurs\bgroup Hartshorne\egroup\egroup{} :
\newblock {\em Algebraic geometry}.
\newblock Graduate Texts in Mathematics, No. 52. Springer-Verlag, New
  York-Heidelberg, 1977.

\bibitem{jan}
J.C. \bgroup\fonteauteurs\bgroup Jantzen\egroup\egroup{} :
\newblock {\em Representations of algebraic groups}.
\newblock Pure and Applied Mathematics, 131. Academic Press, Inc., Boston, MA,
  1987.

\bibitem{kan}
S.~S. \bgroup\fonteauteurs\bgroup Kannan\egroup\egroup{} :
\newblock Projective normality of the wonderful compactification of semisimple
  adjoint groups.
\newblock {\em Math. Z.}, 239(4)\string:\penalty500\relax 673--682, 2002.

\bibitem{lu}
D.~\bgroup\fonteauteurs\bgroup Luna\egroup\egroup{} :
\newblock Toute variété magnifique est sphérique.
\newblock {\em Transform. Groups}, 3(3)\string:\penalty500\relax 249--258,
  1996.

\bibitem{was}
B.~\bgroup\fonteauteurs\bgroup Wasserman\egroup\egroup{} :
\newblock Wonderful varieties of rank two.
\newblock {\em Transform. Groups}, 1\string:\penalty500\relax 375--403, 1996.

\end{thebibliography}

\end{document}